\documentstyle[12pt]{article}

\renewcommand{\theequation}{\thesection.\arabic{equation}}

\begin{document}
\baselineskip=20pt
\renewcommand{\thefootnote}{\fnsymbol{footnote}}
\newpage
\pagestyle{empty}
\setcounter{page}{0}

\vfill
\begin{center}

{\LARGE {\bf {\sf
A nested
sequence of projectors and corresponding braid matrices $\hat
R
(\theta)$:
({\rm 1}) Odd dimensions. }}}
\\[0.8cm] {\large
A.Chakrabarti

{\em

Centre de Physique Th\'eorique\footnote{Laboratoire
Propre
du CNRS UPR A.0014}, Ecole Polytechnique, 91128 Palaiseau Cedex,
France.\\
e-mail
chakra@cpht.polytechnique.fr}}

\end{center}

\smallskip

\smallskip

\smallskip

\smallskip

\smallskip

\smallskip

\begin{abstract}

A basis of
$N^2$ projectors, each an ${N^2}\times{N^2}$ matrix with constant
elements, is implemented to construct a class of braid matrices
$\hat{R}(\theta)$, $\theta$ being the spectral parameter. Only odd
values of $N$ are considered here. Our ansatz for the projectors
$P_{\alpha}$ appearing in the spectral decomposition of
$\hat{R}(\theta)$ leads to exponentials $exp(m_{\alpha}\theta)$ as
the coefficient of $P_{\alpha}$. The sums and differences of such
exponentials on the diagonal and the antidiagonal respectively provide
the $(2N^2 -1)$ nonzero elements of $\hat{R}(\theta)$. One element at
the center is normalized to unity. A class of supplementary constraints
imposed by the braid equation leaves $\frac{1}{2}(N+3)(N-1)$ free
parameters $m_{\alpha}$. The diagonalizer of $\hat{R}(\theta)$ is
presented for all $N$. Transfer matrices $t(\theta)$ and $L(\theta)$
operators corresponding to our $\hat{R}(\theta)$ are studied. Our
diagonalizer signals specific combinations of the components of the
operators that lead to a quadratic algebra of $N^2$ constant $N\times N$
matrices. The $\theta$-dependence factors out for such combinations.
$\hat R(\theta)$ is developed in a power series in $\theta$. The basic
difference arising for even dimensions is made explicit. Some special
features of our $\hat{R}(\theta)$ are discussed in a concluding section.
\end{abstract}

\vfill
\newpage

\pagestyle{plain}

\section{Introduction}

In Sec.$8$ of Ref.$1$ a sequence of projectors with constant
elements and
particularly
simple and convenient properties were introduced
for arbitrary dimension
$N$ ( i,e,
$N^2\times N^2$ matrices ). For the case
$N =2$ they provide the spectral
resolutions of
the $6$-vertex and the
$8$-vertex models (Sec.$6$ and Sec.$7$ of Ref.$1$,
citing other
sources).
Along with $N$ the set of projectors is enlarged in number and
in
dimension
systematically at each step to give what we called a "nested
sequence". The
projectors
were presented Ref.$1$ for all $N$ and their
basic features were studied,
including
diagonalization, for arbitrary $N$.
But no higher dimensional braid
matrices were
constructed on such  bases.
That was "beyond the scope" of that paper. Here
we enlarge the
scope and
present explicit constructions for all $\it odd$ $N$. Even dimensions
will
be studied separately elsewhere. Such a separation corresponds
to
strikingly
different features arising in the respective cases.

 After
obtaining explicitly $\hat{R}(\theta)$ the corresponding
transfer
matrices
$t(\theta)$ and $L(\theta)$ operators are studied. They
are found to lead
to a remarkable
class of quadratic algebra (Sec.$4$).
Development of $\hat R(\theta)$ in
powers of the
spectral parameter
$\theta$ is also studied (Sec.$5$). Basic differences
arising for
even
dimensions are pointed out. (Sec.$6$). The special features of our
class of
solutions are
discussed in  conclusion (Sec.$7$). Construction of
our solutions is
presented in App.$A$
and some basic results concerning
$t(\theta)$ and $L(\theta)$ are collected
together
in
App.$B$.

\section{Braid matrices for odd dimensions (Ansatz and
Solutions) :}
\setcounter{equation}{0}
We start by specifying our notations and conventions in detail sinc
e they
turn out to be
crucial in successful construction of the solutions.
Thus, rather than
using the simple and
elegant notation of Sec.$8$ of
Ref.$1$  for our projectors we intorduce
below a srtucture
better suited to
our present purpose.

 Let
$$N=2p - 1 \qquad ( p=2,3,...)$$
and $$ \bar{i}
=N-i+1 \qquad   (i+ \bar{i}=2p, \quad \bar{\bar{i}}=i)$$

so that for
$$i=1,2,...,(p-1)$$
respectively $$\bar{i}= (2p-1),(2p-2),...,(p+1)$$
and
$$\bar{p}=p$$

The $N^2 \times N^2$ braid matrix $\hat{R}(\theta)$, with
the spectral
parameter $\theta$,
is given in terms of its components
as
\begin{equation}
\hat{R}(\theta) = ({\hat{R}}(\theta))_{ab,cd}
(ab)\otimes (cd)
\end{equation}
where $(a,b,c,d)$ take values in in the
domain $(i,\bar{i},p)$ and $(ab)$
is the $N\times
N$ matrix with only one
nonzero element,unity, at row $a$ and column $b$.

 The basis of projectors
is given by ( with $\epsilon = \pm $) the set

$$ P_{pp}=
(pp)\otimes(pp)$$
$$2 P_{pi(\epsilon)}= (pp)\otimes \biggl( ((ii)
+({\bar{i}} {\bar{i}}))
+\epsilon
((i\bar{i})+(\bar{i}i))\biggr)$$
$$2
P_{ip(\epsilon)}=  \biggl( ((ii) +({\bar{i}} {\bar{i}}))
+\epsilon
((i\bar{i})+(\bar{i}i))\biggr) \otimes
(pp)$$
$$2P_{ij(\epsilon)}= \biggl( ((ii)\otimes (jj)+(\bar{i}
\bar{i})\otimes
(\bar{j}\bar{j}))
+\epsilon ((i\bar{i})\otimes (j\bar{j})+
(\bar{i}i)\otimes
(\bar{j}j))\biggr)$$
\begin{equation}
2P_{i\bar{j}(\epsilon)}= \biggl(
((ii)\otimes (\bar{j}\bar{j})+(\bar{i}
\bar{i})\otimes
(jj)) +\epsilon
((i\bar{i})\otimes (\bar{j}j)+
(\bar{i}i)\otimes
(j\bar{j}))\biggr)
\end{equation}

 Condensing the
triplets $(i,j,\epsilon),(i,p,\epsilon),..$ and also
$(pp)$
into
$(\alpha,\beta,...)$ the basis $(2.2)$
satisfies
\begin{equation}
P_{\alpha}P_{\beta} =
{\delta}_{\alpha\beta}P_{\alpha},\qquad
\sum_{\alpha}P_{\alpha}=I_{{N^2}\times {
N^2}}
\end{equation}
 The total number of $P_{\alpha}$ is
$$
1+4(p-1)+4(p-1)^2 = (2p-1)^2 =N^2$$
They have, apart from the overall
factor $\frac{1}{2}$ for all projectors
except $P_{pp}$,
only the constant
elements $(\pm1,0)$. There is, for example, no $q$ in our
formalism.
The
braid matrix is $\it{postulated}$ in the spectrally resolved
form
\begin{eqnarray}
\nonumber
\hat{R}(\theta) & = P_{pp} +
\sum
_{i,\epsilon}\bigl(f^{(\epsilon)}_{pi}(\theta)
P_{pi(\epsilon)} +
f^{(\epsilon)}_{ip}(\theta)P_{ip(\epsilon)}\bigr)  \\
 & + \sum
_{i,j,\epsilon} \bigl( f^{(\epsilon)}_{ij}(\theta) P_{ij(\epsilon)}
+
f^{(\epsilon)}_{i\bar{j}}(\theta)P_{i\bar{j}(\epsilon)}\bigr)  \\
\nonumber
\end{eqnarray}

The coefficient of $P_{pp}$ is normalized to
unity. The $(N^2 -1)$ functions
$f^{(\epsilon)}_{ab}$ are to be extracted
from the constraints imposed by
the
braid equation
\begin{equation}
\hat{R}_{12}(\theta)\hat{R}_{23}(\theta+{\theta}')\hat{R}_{12}({\theta}')=
\hat{R}_{23}({\theta}')\hat{R}_{12}({\theta+\theta}')
\hat{R}_{23}(\theta)
\end{equation}
Here (suppressing $\theta$)

$$\hat{R}_{12} = \hat{R} \otimes I_{N\times N}, \qquad
\hat{R}_{23} =I_{N\times N}\otimes \hat{R}$$

In terms of the coefficients
$({\hat{R}}(\theta))_{ab,cd}$ defined in
$(2.1)$ one obtains
( summing over
the repeated indices $(l,m,n)$
)
\begin{eqnarray}
\nonumber
(\hat{R}(\theta))_{al,cm}(\hat{R}(\theta+\theta'))_{mn,ef}
(\hat{R}(\theta'))_{lb,nd} \\
=(\hat{R}(\theta'))_{cl,em}(\hat{R}(\theta+\theta'))_{ab,ln}
(\hat{R}(\theta))_{nd,mf}
\end{eqnarray}

This corresponds to the point $(ab)\otimes
(cd)\otimes(ef)$ of the base
space $V\otimes V
\otimes V$. Our ansatz
$(2.4)$ along with $(2.2)$ implies very srtong
constraints (typical
of $\it
odd$ dimensions). The solutions are obtained in App.$A$. One
has
\begin{equation}
 f^{(\epsilon)}_{ab}(\theta) = exp
(m^{(\epsilon)}_{ab} \theta) ; \qquad (ab)
=
(pi),(ip),(ij),(i\bar{j})
\end{equation}

where the parameters
$m^{(\epsilon)}_{ab}$ are all independent {\it except}
that for
each
$i$
\begin{equation}
m^{(\epsilon)}_{ij}=m^{(\epsilon)}_{i\bar{j}},
\qquad (\bar{j} = 2p - j)
\end{equation}

The constraints $(2.7),(2.8)$ are
{\it necessary and sufficient}. Thus for
$N=3$ one
has

\begin{equation}
\hat{R}  (\theta)  =  \pmatrix{
   a_{+} &0 &0 &0&0
&0 &0&0&a_{-} \cr
   0 &b_{+} &0 &0&0 &0 &0&b_{-}&0 \cr
   0 &0 &a_{+} &0&0
&0 &a_{-}&0&0 \cr
   0 &0 &0 &c_{+}&0 &c_{-} &0&0&0 \cr
   0 &0 &0 &0&1 &0
&0&0&0 \cr
   0 &0 &0 &c_{-}&0 &c_{+} &0&0&0 \cr
   0 &0 &a_{-} &0&0 &0
&a_{+}&0&0 \cr
   0 &b_{-} &0 &0&0 &0 &0&b_{+}&0 \cr
    a_{-}&0 &0 &0&0 &0
&0&0& a_{+} \cr
}
\end{equation}

where
\begin{eqnarray}
\nonumber
a_{\pm}=\frac{1}{2} (
e^{m_{11}^{(+)}{\theta}} \pm e^{m_{11}^{(-)}{\theta}}
)
\\
\nonumber
b_{\pm}=\frac{1}{2} ( e^{m_{12}^{(+)}{\theta}} \pm
e^{m_{12}^{(-)}{\theta}}
) \\
c_{\pm}=\frac{1}{2} (
e^{m_{21}^{(+)}{\theta}} \pm e^{m_{21}^{(-)}{\theta}} )
\end{eqnarray}

The
six parameters remaining after application of $(2.8)$ (which imposes
the
repetation
of $a_{\pm}$) are all ${\it independent}$. For $
m_{ab}^{(-)}=
m_{ab}^{(+)}$ one obtains
hyperbolic functions as particular
cases. For all $N$, the nonzero elements
are confined
to {\it the diagonal
and the antidiagonal} as above with a common element,
unity, at the
centre.
Apart from the normalized element, the coefficients of the
projectors are
simply
exponentials. The total number of independent parameters
$m_{ab}^{\pm}$ is

\begin{equation}
(2p -1)^2 -1 - 2(p-1)^2 \quad =2(p^2
-1) \quad = \frac{1}{2}(N+3)(N-1)
\end{equation}

Note that the coefficient
of $P_{pp}$ in $(2.4)$ has to be nonzero for
$\hat{R}(\theta)$
to be
invertible and hence can safely be normalized to unity. Indeed,
each
coefficient in
$(2.4)$ has to be nonzero for $\hat R(\theta)$ to be
invertible. This is
more evident after
diagonalization (Sec.$3$). For even
$N$ there is   ${\it no}$ index $p=
\bar p$. In
App.$A$ the crucial role of
the index $p$ will be made more evident. The
projectors
$P_{pi(\epsilon)}$
and $P_{ip(\epsilon)}$ will be seen to impose the
highly
constrained
solutions $(2.7)$ with $(2.8)$.

Implementing $(2.8)$
one obtains from $(2.4)$
\begin{equation}
\hat{R}(\theta) = P_{pp} + \sum
_{i,\epsilon}\bigl( f^{(\epsilon)}_{pi}(\theta)
P_{pi(\epsilon)} +
f^{(\epsilon)}_{ip}(\theta)P_{ip(\epsilon)}\bigr)+\sum
_{i,j,\epsilon}f^{(\epsilon)}_{ij}(\theta)\bigl(P_{ij(\epsilon)}+
P_{i\bar{j}(\epsilon)}\bigr)
\end{equation}

Defining
\begin{equation}
\tilde{
P}_{ij(\epsilon)}= P_{ij(\epsilon)} +
P_{i\bar{j}(\epsilon)}
\end{equation}
and conserving all other projectors
as before one obtains a basis of $(2p^2
-1)$
projectors still satisfying
$(2.3)$ where now the indices summed over
are
$(i,j,\epsilon),(i,p,\epsilon),(p,i,\epsilon),(pp)$.

Now

\begin{equation}
\hat{R}(\theta) = P_{pp} + \sum _{i,\epsilon}\bigl(
f^{(\epsilon)}_{pi}(\theta)
P_{pi(\epsilon)} +
f^{(\epsilon)}_{ip}(\theta)P_{ip(\epsilon)}\bigr)+\sum
_{i,j,\epsilon}
f^{(\epsilon)}_{ij}(\theta)\tilde{P}_{ij(\epsilon)}
\end{equation}

In this basis {\it all} the $\frac{1}{2}(N+3)(N-1)$ parameters are
independent. When they
are all chosen to be distinct (and different
from $1$) the polynomial
equation ( of
$\frac{1}{2}(N+3)(N-1)$ degree and
with distinct roots ) satisfied by $\hat
{R}(\theta)$
and the projectors in
terms of $\hat {R}(\theta)$ are obtained respectively
as in
$(1.5)$and
$(1.6)$ of Ref.$1$. The initial basis, due to the symmetry and
simplicity
of
the projectors, is most convenient for certain purposes. The second one
has
the virtue of
eliminating constraints. Each should be implemented
according to the context.

 If two or more of the $\frac{1}{2}(N+3)(N-1)$
free parameters are allowed
to coincide,
then introducing the sum of the
corresponding projectors (as in $(2.13)$)
the basis can
again be redefined
(as in $(2.14)$). The degree of the minimal polynomial
equation
satisfied
by $\hat R(\theta)$ diminishes correspondingly.

 Our matrices all
satisfy

\begin{equation}
\hat R(-\theta)\hat R(\theta) =I, \qquad \hat
R(0)
=I
\end{equation}

\section{Diagonalization:}
\setcounter{equation}{0}

 Our general approach to
diagonalization is presented step by step in Sec.$9$ of Ref.$1$.

 The
matrix $M$ that diagonalizes each projector $P_{\alpha}$ of $(2.3)$
(namely, $P_{pp}$,
$P_{pi(\epsilon)}$, $ P_{ip(\epsilon)}$, $
P_{ij(\epsilon)}$,
$ P_{i{\bar j}(\epsilon)}$ of $(2.2)$ ) and hence
$\hat{R}(\theta)$ of
$(2.4)$ is
given below. As compared to the the results
of Sec.$8$ of Ref.$1$, $M$ is
presented
here in our current notations.

Set
  $$\sqrt{2} M = \sqrt{2} M^{-1}=\sqrt{2} (pp)\otimes (pp) + $$
$$
(pp)\otimes \biggl ( \sum_{i}
\bigl( (ii)-({\bar i}{\bar i})+(i {\bar i}) +
({\bar i}i) \bigr)\biggr)
 + \biggl ( \sum_{i}\bigl( (ii)-({\bar i}{\bar
i}) + (i {\bar i})+({\bar
i}i) \bigr)\biggr)
\otimes (pp) +
$$
\begin{equation}
  \sum _{i,j }\biggl ( \bigl( (ii)-({\bar i}{\bar
i})\bigr)\otimes \bigl(
(jj)+({\bar
j}{\bar j})\bigr) +\bigl((i {\bar i}) +
({\bar i}i) \bigr) \otimes \bigl(
(j{\bar j}) +
({\bar j}j)\bigr)
\biggr)
\end{equation}

One verifies in a straightforward fashion (with
$\epsilon =\pm 1$ on the
right) that

$$ M P_{pp}M^{-1} = (pp)\otimes (pp)
$$

$$ 2 M P_{pk(\epsilon)}M^{-1} =
(pp)\otimes
\bigl((1+\epsilon)(kk)+(1-\epsilon)(\bar{k}\bar{k}) \bigr)$$
$$
2 M P_{kp(\epsilon)}M^{-1}
=
\bigl((1+\epsilon)(kk)+(1-\epsilon)(\bar{k}\bar{k}) \bigr) \otimes
(pp)$$
$$ 2 M P_{kl(\epsilon)}M^{-1} =
(1+\epsilon)(kk)\otimes
(ll)+(1-\epsilon)(\bar{k}\bar{k})\otimes
(\bar{l}\bar{l})
$$
\begin{equation}
 2 M P_{k\bar{l}(\epsilon)}M^{-1}
=
(1+\epsilon)(kk)\otimes
(\bar{l}
\bar{l})+(1-\epsilon)(\bar{k}\bar{k})\otimes (ll)
\end{equation}

Hence taking account of $(2.8)$ (i,e,
$f^{(\epsilon)}_{ij}(\theta)
=f^{(\epsilon)}_{i \bar{j}}(\theta)$) one
obtains

$$ 2 M \hat {R} (\theta) M^{-1}  = 2(pp)\otimes (pp) $$
$$+\sum
_{i,\epsilon}\biggl(f^{(\epsilon)}_{pi}(\theta)
((1+\epsilon)(pp)\otimes
(ii)
+(1-\epsilon)(pp)\otimes(\bar{i}\bar{i}))$$
$$+f^{(\epsilon)}_{ip}(\theta)
((1+\epsilon)(ii)\otimes (pp)
+(1-\epsilon)(\bar{i}\bar{i})\otimes (pp))
\biggr) $$
$$\sum_{i,j,\epsilon}
\biggl (f^{(\epsilon)}_{ij}(\theta)
\bigl((1+\epsilon)((ii)\otimes (jj)
+(ii)\otimes (\bar{j}\bar{j}))
$$
\begin{equation}
+(1-\epsilon)((\bar{i}\bar{i})\otimes
(jj)
+(\bar{i}\bar{i})\otimes(\bar{j}\bar{j})))\biggr)
\end{equation}

 For
$N=3$ this gives
$$ M \hat {R} (\theta) M^{-1} \equiv \hat {R}_{d}
(\theta)$$
\begin{equation}
= ( e^{m_{11}^{(+)}\theta},
e^{m_{12}^{(+)}\theta}, e^{m_{11}^{(+)}\theta},
e^{m_{21}^{(+)}\theta},1,
e^{m_{21}^{(-)}\theta},
e^{m_{11}^{(-)}\theta},
e^{m_{12}^{(-)}\theta},
e^{m_{11}^{(-)}\theta})_{(diag)}
\end{equation}

The diagonalizer is

\begin{equation}
\sqrt 2 M = \sqrt 2 M^{-1} =  \pmatrix{
   1
&0 &0 &0&0 &0 &0&0&1 \cr
   0 &1 &0 &0&0 &0 &0&1&0 \cr
   0 &0 &1 &0&0 &0
&1&0&0 \cr
   0 &0 &0 &1&0 &1 &0&0&0 \cr
   0 &0 &0 &0&\sqrt 2 &0 &0&0&0
\cr
   0 &0 &0 &1&0 &-1 &0&0&0 \cr
   0 &0 &1 &0&0 &0 &-1&0&0 \cr
   0 &1
&0 &0&0 &0 &0&-1&0 \cr
   1 &0 &0 &0&0 &0 &0&0&-1 \cr
}
\end{equation}

 The generalizations of $(3.4)$ and $(3.5)$ for all $N$
are quite evident.

 If an $\hat R(\theta)$ satisfying the braid equation
$(2.5)$ is
diagonalized the
corresponding $\hat R_{d}(\theta)$, in general,
does ${\it not}$ directly
satisfy $(2.5)$.
This is evident from all the
examples of Ref.$1$. The general explanation is
simple. Interpolated
factors of the type $M_{12}M^{-1}_{23}$ will be
lacking in the latter case
as compared to the former. If $\hat {R} (\theta)$ is diagonal {\it to
start with } $(2.6)$ reduces to
$$(\hat{R}(\theta))_{aa,bb}(\hat{R}(\theta+{\theta}'))_{bb,cc}
(\hat{R}({\theta}'))_{aa,bb}$$
\begin{equation}
=(\hat{R}({\theta}'))_{bb,cc}(\hat{R}({\theta+\theta}'))_{aa,bb}
(\hat{R}(\theta))_{bb,cc}
\end{equation}
The braid equation is satisfied if, for each
$(a,b)$,

\begin{equation}
(\hat{R}(\theta))_{aa,bb}(\hat{R}({\theta}'))_{aa,bb}=
(\hat{R}({\theta+\theta}'))_{aa,bb}
\end{equation}
i,e, if
\begin{equation}
(\hat{R}(\theta))_{aa,bb}=
e^{m_{ab}\theta}
\end{equation}
where the parameters $m_{ab}$ are mutually
${\it independent}$. Now,
conversely, if
$\hat{R}(\theta)$ is conjugated
as
\begin{equation}
{\hat{R}}'(\theta)
=A\hat{R}(\theta)A^{-1}
\end{equation}
in general, ${\hat{R}}'(\theta)$
will  no longer satisfy the braid equation
since
such products as
$A^{-1}_{12}A_{23}$ will depend on the srtucture of $A$.
The srtucture
of
our
$M$ is such that for arbitrary odd $N$
$$ M^{-1} \hat {R}_{d}
(\theta) M $$
 continues to satisfy the braid equation ${\it provided
}$
$$m^{(\epsilon)}_{ij}=m^{(\epsilon)}_{i\bar{j}}$$.

 Thus it is seen how
the $2(p-1)^2$ crucial constraints $(2.8)$, the
structure of our
nested
sequence of projectors and that of our $M$ are all linked.

 The relevance
of our $M$ to the algebra of the $L$-operators is pointed
out at the end
of
Sec.$4$ after displaying the crucial algebraic srtucture arising
there.

\section{$L(\theta)$-operators and transfer
matrices:}
\setcounter{equation}{0}

  A general discussion, citing relevant sources, is presented in
App.$B$.
Here the basic
results concerning the the $N\times N$
realizations of the $N^2$ blocks of
the transfer
matrix $t(\theta)$ and the
operator $L^{+}(\theta)$ are used in the context
of  braid
matrices
constructed in Sec.$2$ and App.$A$.

 In $(B.22)$ and $(B.23)$ we show in a
transparent fashion why, unless
$(B.16)$ is
generalized, say, by
implementing central operators in the argument of
$\hat R (\theta
-
\theta')$, one cannot obtain an $L^{-}(\theta) \neq L^{+}(\theta)$. We
do
not study such
general structures here and hence consider only the
above-mentioned fundamental
realizations of $L^{+}(\theta)$ with the
standard prescription for
coproduct. This will,
in any case provide a
subalgebra in an appropriately generalized quasi-Hopf
structure.
This
$L^{+}(\theta)$ and
$t(\theta)$, as shown in $(B.28)$, are related (
for the fundamental
$N\times N$
representations of blocks ) as
\begin{equation}
t(\theta) =P L^{+}(\theta) P, \qquad
(t_{ab}(\theta))_{cd}=(L^{+}_{cd}(\theta))_{ab}
\end{equation}
In studying multistate statistical models corresponding to our
$\hat R(\theta)$ ( see
the comments and references in Sec.$7$ ) the algebra of the blocks of
$t(\theta)$ is
particularly relevant. In our case this
algebra is found ( see below ) to be very simply related to the
corresponding one for $L^{+}(\theta)$. So one can start either with
$L^{+}(\theta)$ or $t(\theta)$ and then obtain the other easily. We
choose to display the remarkable structure that emerges first in terms of
$L^{+}(\theta)$. We start with $(B.14)$, i,e,
\begin{equation}
L^{+}(\theta)=\hat R(\theta) P
\end{equation}
In terms of the matrices $(ab)$ defined below $(2.1)$, one obtains
$$L^+_{(pp)}=(pp) \equiv X_{pp}$$
$$e^{-m_{ip}^{(\epsilon)}\theta}(L_{ip}^{+}{(\theta)}+{\epsilon}
L_{\bar{i}p}^{+}{(\theta)})= (pi) + \epsilon (p\bar i) \equiv
X_{pi}^{(\epsilon)}$$
$$e^{-m_{pi}^{(\epsilon)}\theta}(
L_{pi}^{+}{(\theta)}+{\epsilon}
L_{p\bar{i}}^{+}{(\theta)})= (ip) + \epsilon (\bar {i}p)
\equiv X_{ip}^{(\epsilon)}$$
$$e^{-m_{ij}^{(\epsilon)}\theta}(
L_{ij}^{+}{(\theta)}+{\epsilon}
L_{\bar{i}\bar{j}}^{+}{(\theta)})= (ij) +
\epsilon (\bar {i}\bar{j})
\equiv
X_{ij}^{(\epsilon)}$$
\begin{equation}
e^{-m_{ij}^{(\epsilon)}\theta}( L_
{i\bar{j}}^{+}{(\theta)}+{\epsilon}
L_{\bar{i}j}^{+}{(\theta)})= (\bar{j}i)
+ \epsilon (j\bar {i}) \equiv
X_{\bar{j}i}^{(\epsilon)}
\end{equation}

 In
the last equation $(2.8)$ has been implemented i,e,
$$
m_{i\bar{j}}^{(\epsilon)} =m_{ij}^{(\epsilon)} $$

 From these one
obtains

$$ 2L_{ij}^{+}(\theta) =
(e^{m_{ij}^{(+)}\theta}X_{ij}^{(+)}+
e^{m_{ij}^{(-)}\theta}X_{ij}^{(-)})$$
\begin{equation}
2L_{\bar{i}\bar{j}}^{+}(\theta) =
(e^{m_{ij}^{(+)}\theta}X_{ij}^{(+)}-
e^{m_{ij}^{(-)}\theta}X_{ij}^{(-)})
\end{equation}
and so on.

For $N=3$ one obtains (with $\epsilon=\pm 1$ in the
matrices on the right)

$$ L_{22}^{+}(\theta) = \pmatrix{
0&0&0
\cr
0&1&0\cr
0&0&0 \cr }$$

$$ L_{12}^{+}(\theta) +\epsilon
L_{\bar{1}2}^{+}(\theta)=
e^{ m_{12}^{(\epsilon)}\theta}\pmatrix{ 0&0&0
\cr
1&0&\epsilon \cr
0&0&0 \cr }$$
$$ L_{21}^{+}(\theta) +\epsilon
L_{2\bar{1}}^{+}(\theta)=
e^{ m_{21}^{(\epsilon)}\theta}\pmatrix{ 0&1&0
\cr
0&0&0 \cr
0&\epsilon&0 \cr }$$
$$ L_{11}^{+}(\theta) +\epsilon
L_{\bar{1}\bar{1}}^{+}(\theta)=
e^{ m_{11}^{(\epsilon)}\theta}\pmatrix{ 1&0&0 \cr
0&0&0 \cr
0&0&\epsilon \cr }$$
\begin{equation}
 L_{1\bar{1}}^{+}(\theta
) +\epsilon  L_{\bar{1}1}^{+}(\theta)=
e^{
m_{11}^{(\epsilon)}\theta}\pmatrix{ 0&0&\epsilon \cr
0&0&0\cr
1&0&0 \cr
}
\end{equation}

 The constant $N\times N$ matrices $
X_{ab}^{(\epsilon)}$, where $X_{pp}$
has only one
nonzero element,unity,
and all the others only two ($(1,1)$ or $(1,-1)$)
specify a
quadratic
algebra. We give below only the {\it nonzero} bilinear products,
all
others
vanishing. Further results, such as commutartors, can be
systematically
obtained from
those below:

$$X_{pp}X_{pp}=X_{pp},\qquad
X_{pp}X_{pi}^{(\epsilon)}=X_{pi}^{(\epsilon)},
\qquad
X_{ip}^{\epsilon}X_{pp}=X_{ip}^{\epsilon}$$
$$X_{pi}^{(\epsilon)}X_{ip}^{
{(\epsilon}')} =
(1+{\epsilon}{\epsilon}')X_{pp},
\qquad
X_{ip}^{(\epsilon)}X_{pj}^{{(\epsilon}')}=
X_{ij}^{(\epsilon{\epsilon}')}
+{\epsilon}X_{\bar{i}j}^{(\epsilon{\epsilon}')}$$
$$X_{pi}^{(\epsilon)}X_{ij}^{{
(\epsilon}')}=
X_{pj}^{(\epsilon{\epsilon}')},
\qquad
X_{ij}^{(\epsilon)}X_{jp}^{{(\epsilon}')}=
X_{ip}^{(\epsilon{\epsilon}')}$$
$$X_{pi}^{(\epsilon)}X_{\bar{i}j}^{{(\epsilon}'
)}=\epsilon
X_{pj}^{(\epsilon{\epsilon}')},
\qquad
X_{\bar{j}i}^{(\epsilon)}X_{ip}^{{(\epsilon}')}=
\epsilon{\epsilon}'X_{jp}^{(\epsilon{\epsilon}')}$$
$$X_{ij}^{(\epsilon)}X_{jk}^{{(\epsilon}')}=
X_{ik}^{(\epsilon{\epsilon}')},
\qquad
X_{ij}^{(\epsilon)}X_{\bar{j}k}^{{(\epsilon}')}=
\epsilon
X_{\bar{i}k}^{(\epsilon{\epsilon}')}$$
\begin{equation}
X_{\bar{j}i}^{(\epsilon)
}X_{ik}^{{(\epsilon}')}=
X_{\bar{j}k}^{(\epsilon{\epsilon}')},
\qquad
X_{\bar{j}i}^{(\epsilon)}X_{\bar{i}k}^{{(\epsilon}')}=
\epsilon
X_{jk}^{(\epsilon{\epsilon}')}
\end{equation}

(No sum over repeated
indices.)

  Note that
\begin{equation}
C_{1} \equiv \frac {1}{N}( X_{pp}+
\sum_{i}  X_{ii}^{(+)} )= \frac {1}{N}
I_{N\times N}
\end{equation}

Hence
$(X_{pp} -C_{1})$ and $ (X_{ii}^{(+)} - 2C_{1})$, along with the
others
give an
algebra of $N^2$ traceless matrices.

  Higher dimensional
realizations are given by the coproducts
\begin{equation}
\Delta L =
L\dot{\otimes} L
\end{equation}
Here $\dot{\otimes}$ implies tensor product
combined with matrix
multiplication. The
prescription can be implemented
repeatedly in a straightforward fashion.
But it leads, in
general, to
reducible structures. A systematic study of extraction
of
irreducible
components is beyond the scope of this paper. Let us,
however, take a
closer look at the
structure of the algebra $(4.6)$ and the
special role of the index $p$.

 The generators without $p$ ( i,e,
$X_{ij}^{(\epsilon)}$,
$X_{\bar{i}k}^{(\epsilon)}$)
form a closed
subalgebra. The generators with a single $p$  (
i,e,
$X_{pi}^{(\epsilon)}$,
$X_{ip}^{(\epsilon)}$) provide a semidirect
product structure with the
preceding set. But
now to close it one has to
extend the first set to a direct product
structure by including
$X_{pp}$.

>From $(4.1)$ and $(4.3)$ it can be shown that $t(\theta)$
and
$L^{+}(\theta)$ are
essentially related through the interchange of the
roles of
$X_{pi}^{(\epsilon)}$ and
$X_{ip}^{(\epsilon)}$. Thus for $N=3$
there is an interchange of $b_{\pm}$
and $c_{\pm}$.
 One  obtains for this
case

  \begin{equation}
  t(\theta)  =  \pmatrix{
   a_{+} &0 &0 &0&0 &0
&0&0&a_{-} \cr
   0 &0&0&c_{+} &0 &c_{-}&0&0&0 \cr
   0 &0 &a_{-} &0&0 &0
&a_{+}&0&0 \cr
   0 &b_{+} &0 &0&0&0&0 &b_{-} &0 \cr
   0 &0 &0 &0&1 &0
&0&0&0 \cr
   0 &b_{-} &0 &0&0&0&0 &b_{+} &0 \cr
   0 &0 &a_{+} &0&0 &0
&a_{-}&0&0 \cr
   0 &0&0&c_{-} &0 &c_{+}&0&0&0 \cr
    a_{-}&0 &0 &0&0 &0
&0&0& a_{+} \cr         }
\end{equation}

\begin{equation}
L^{+}(\theta)  =  \pmatrix{
   a_{+} &0 &0 &0&0 &0 &0&0&a_{-} \cr
   0
&0&0&b_{+} &0 &b_{-}&0&0&0 \cr
   0 &0 &a_{-} &0&0 &0 &a_{+}&0&0 \cr
   0
&c_{+} &0 &0&0&0&0 &c_{-} &0 \cr
   0 &0 &0 &0&1 &0 &0&0&0 \cr
   0 &c_{-}
&0 &0&0&0&0 &c_{+} &0 \cr
   0 &0 &a_{+} &0&0 &0 &a_{-}&0&0 \cr
   0
&0&0&b_{-} &0 &b_{+}&0&0&0 \cr
    a_{-}&0 &0 &0&0 &0 &0&0& a_{+} \cr
}
\end{equation}
 From these the $3\times 3$ blocks can be read off.

   We
close this section by pointing out the relevance of our diagonalizer
$M$ of
Sec.$3$
to the structure of $ L^{+}(\theta)$ ( and hence of $t(\theta)$ ).
If one
constructs
\begin{equation}
ML_{+}(\theta) M^{-1}
\end{equation}

precisely the combinations on the left of the set $(4.3)$ are seen
to
emerge. Thus our
$M$ leads directly to the remarkable structure
$(4.6)$.

\section{$\theta$-expansion:}
\setcounter{equation}{0}

 Let us start with the following
notations and conventions:

 $(1)$: The condensed notation
$(\alpha,\beta,...)$ of $(2.3)$ implies for
each index
$\alpha$
either
$(pp)$ or a triplet $(p,i,\epsilon), (i,j,\epsilon),...$ and so on.
We
introduce sum over
${\alpha}'$ where in $\sum_{\alpha'}$ the index
$(pp)$ is $\it excluded$.
As for the
other projectors one may consider
alternatively either the basis given by
$(2.12)$ or
that by
$(2.14)$.

$(2)$: We also define
\begin{equation}
H \equiv \sum_{\alpha'}
m_{\alpha'} P_{\alpha'}
\end{equation}

when ,using
$(2.3)$,
\begin{equation}
H^n = (\sum_{\alpha'} m_{\alpha'} P_{\alpha'})^n
= \sum_{\alpha'}
m^n_{\alpha'} P_{\alpha'}
\end{equation}

Now one can
expand as follows ( with $n \geq 1$ )

$$\hat R(\theta) = P_{pp}
+\sum_{\alpha'} e^{m_{\alpha'}\theta} P_{\alpha'}$$
$$ = P_{pp}
+\sum_{\alpha'}\biggl(1+\sum_n
\frac{(m_{\alpha'}\theta)^n}{n!}\biggr)
P_{\alpha
'}
\qquad = I +\sum_{n} \frac{{\theta}^n}{n!} \bigl(
\sum_{\alpha'}
m_{\alpha'}^{n}
P_{\alpha'}\bigr)$$
\begin{equation}
=I +
\sum_{n} \frac{{\theta}^n}{n!}H^n  \qquad =e^{\theta H}
\end{equation}

Addition of $P_{pp} =(I - \sum_{\alpha'} P_{\alpha'}) $ to $H$
corresponds
to a change of
normalization of $\hat R(\theta)$ along with an
evident redefinition
$m_{\alpha'} \rightarrow (m_{\alpha'} -1)$. None of
the considerations
below are affected
by such a redefinition
($\sum_{\alpha'} \rightarrow  \sum_{\alpha}$) of $H$.
  More generally, say
for $q$- deformed $(A,B,C,D)$-type algebras, if $\hat
R (\theta)$
is
spectrally resolved  on a complete basis of projectors (Sec.$2$,
Ref.$1$),
setting $l_i(\theta) =ln k_i (\theta)$ and normalizing suitably
one
obtains, following the
steps leading to $(5.3)$,
\begin{equation}
\hat
R(\theta) = \sum_{i} k_{i}(\theta) P_i \qquad = \sum_{i}e
^{l_{i}(\theta)}
P_i \qquad
= e^{(\sum_{i} l_{i}(\theta) P_i)}
\end{equation}
 Here, in
general, upon expansion in powers of $\theta$ the exponents
$l_i(\theta)$
lead
to fairly involed structures. In our present case $\theta$ is simply
a
factor in the
exponent. Hence the situation is much simpler.
 Using
$(5.3)$ the braid equation becomes (with $H_{12} = H\otimes I$ ,
$H_{23}=
I\otimes
H$)
\begin{equation}
e^{\theta H_{12}}e^{(\theta +\theta')
H_{23}}e^{\theta' H_{12}} =
 e^{\theta ' H_{23}}  e^{(\theta +\theta')
H_{12}}e^{\theta H_{23}}
\end{equation}
Setting, with $(n,n',n'') \geq
1$,
\begin{equation}
S = \sum_{n} \frac{{\theta}^n}{n!}H_{12}^n ,\qquad
S'=\sum_{n'}
\frac{{\theta'}^{n'}}{n'!}H_{12}^{n'},
\qquad S'' = \sum_{n''}
\frac{{(\theta +\theta')}^{n''}}{n''!}H_{23}^{n''}
\end{equation}

The left
hand side of $(5.5)$ is
$$(L.H.S.) = (I+S)(I+S'')(I+S')
$$
\begin{equation}
= I + (S+S'+S'')+ ( SS'+ SS''+ S''S') +
SS''S'
\end{equation}
The $(R.H.S.)$ is obtained from the $(L.H.S.)$ via
the following interchanges :
\begin{equation}
  (12)  \leftrightarrow (23),
\qquad \theta \leftrightarrow \theta'
\end{equation}
 Now let us compare
the coefficients of $\theta ^r \theta'^s$ for different
pairs $(r,s)$
on
both sides of $(5.5)$.

  The linear and the quadratic terms on both sides
are found to be
symmetric under $(5.8)$
and hence cancel. Among the cubic
terms only the coefficients of
            $$ \theta \theta' (\theta
+\theta')$$
are found to lead to a nontrivial relation. One obtains, on
regrouping terms,
\begin{equation}
[[H_{12},H_{23}],H_{12}] =
[[H_{23},H_{12}],H_{23}]
\end{equation}
Compare this with $(2.5)$. See also
the remarks in Sec.$7$.

 But from $(5.2)$ one
obtains
\begin{equation}
H_{12}^2 = \sum_{\alpha'}m_{\alpha'}^2
(P_{\alpha'})_{12}, \qquad H_{23}^2 =
\sum_{\alpha'}m_{\alpha'}^2
(P_{\alpha'})_{23}
\end{equation}
 Hence in terms of the projectors one
obtains
$$\sum_{\alpha',\beta',\gamma'}m_{\alpha'}m_{\beta'}m_{\gamma'}
\biggl(
(P_{\alpha'})_{12}(P_{\beta'})_{23}(P_{\gamma'})_{12}
-(P_{\alpha'})_{23
}(P_{\beta'})_{12}(P_{\gamma'})_{23} \biggr)$$
\begin{equation}
 =
\frac{1}{2}\sum_{\alpha',\beta'}m_{\alpha'}m_{\beta'}(m_{\alpha'}
-
m_{\beta'})
\biggl((P_{\alpha'})_{12}(P_{\beta'})_{23}
-(P_{\alpha'})_{23}(P_{
\beta'})_{12} \biggr)
\end{equation}
Since there are
$\frac{1}{2}(N+3)(N-1)$ independent parameters
$m_{\alpha}$,
comparing
coefficients of distinct triplets on each side one obtains a
series of
results. We will
not display them explicitly.
 In Sec.$3$ of
Ref.$2$ we have studied analogous reductions ( from
trilinear to
bilinear
forms ) for $q$-deformed unitary, orthogonal and symplectic cases.
There they
were studied in the context of "modified braid equations"
($[3],[4]$)
presented as a complementery facet of Baxterization (i,e, the
introduction
of a spectral
parameter ). Here we {\it started} from the
$\theta$-dependent form $(2.5)$ and
implemented our $\theta$- expansion
leading to the hierarchy starting
with$(5.9)$ and
$(5.11)$. Without
attempting to analyse how the higher order members of the
hierarchy can be
reduced in order, in successive steps, we just mention the
following
point
concerning $(5.11)$.

 In $(1.18)$ of Ref.$2$ , even for the orthogonal and
the symplectic cases the
modified braid equation could be expressed in
terms of tensored $(\hat
R(\theta))^{\pm 1}$ by expressing the projectors
in their terms using the
minimal (cubic)
polynomial equation satisfied by
$\hat R(\theta)$ . For the unitary case (
with a
quadratic polynomial ) the
task was much more simple. In our present case,
despite
various
particularly simple aspects, the order of the minimal polynomial
increases
as
$N^2$ instead of remaining fixed as for the cases mentioned before.
Hence
relations of
the type $(5.11)$ are best considered in terms of
projectors themselves.

     Expansions in terms of the spectral parameter
has been  considered in
the context
of Yangian Double and central
extensions $[5,6]$. We intend to study elsewhere
analogous aspects
generalizing our class of braid matrices.

\section{Comparison with even
dimensional
cases:}
\setcounter{equation}{0}

 The sequence of projectors presented in Sec.$8$  of Ref.$1$ is a
direct
generalization of
the basis arising in the spectral resolution of
the $6$-vertex and the
$8$-vertex braid
matrices. From the abundant
literature on such models the most directly
relevant sources
are cited in
Sec.$6$ and Sec.$7$ of Ref.$1$. These $4\times 4$ projectors
are (
with
$\epsilon = \pm 1$ in the matrices )

\begin{equation}
 2
P_{1(\epsilon)} =   \pmatrix{
   1 &0 &0 & \epsilon \cr
   0&0 &0 &0 \cr
0 &0 &0&0 \cr\epsilon &0&0&1      }, \qquad  2 P_{2(\epsilon)} =
\pmatrix{
   0 &0 &0 & 0\cr
   0&1 &\epsilon &0 \cr
   0 &\epsilon &1&0
\cr
   0&0&0&0     }
\end{equation}

  But even for this simplest member $(
N=2 )$ of the hierarchy the
coefficients
in
\begin{equation}
\hat{R}(\theta) =  \sum _{\epsilon}\bigl(
f_{1(\epsilon)}(\theta)
P_{1(\epsilon)} +
f_{2(\epsilon)}(\theta)P_{2(\epsilon)}\bigr)
\end{equation}
 are not
constrained to simple exponentials as for $N=(2p-1)$. For the
$8$-vertex
model
( see sources cited in Sec.$7$ of Ref.$1$ ) one
obtains
\begin{equation}
f_{1(\pm)} (\theta) = \frac{g_{(\pm)}
(\theta)}{g_{(\pm)} (-\theta)} \qquad
f_{2(\pm)}
(\theta) = \frac{h_{(\pm)}
(\theta)}{h_{(\pm)} (-\theta)}
\end{equation}

where with $z=e^{\theta}$,
two parameters $p$ and $q$ and
\begin{equation}
(x;a)_{\infty} = \prod _{n
\geq 0} (1- xa^n)
\end{equation}
\begin{equation}
g_{\pm}(z) = (\mp
p^{\frac {1}{2}}q^{-1}z;p)_{\infty}(\mp
p^{\frac
{1}{2}}qz^{-1};p)_{\infty}
\end{equation}
\begin{equation}
h_{\pm}(z) =
( q^{\frac {1}{2}} z^{-\frac {1}{2}}\pm q^{-\frac {1}{2}}
z^{\frac
{1}{2}})
(\mp pq^{-1}z;p)_{\infty}(\mp p
qz^{-1};p)_{\infty}
\end{equation}
 The question of normalization is
discussed in Sec.$7$ of Ref.$1$.
 In the trigonometric $6$-vertex limit one
obtains(as in Sec.$6$ of Ref.$1$)
\begin{equation}
f_{1(\pm)}(\theta) =1,
\quad f_{2(+)}(\theta) = \frac {cosh \frac {1}{2}(
\gamma -
\theta)}{cosh
\frac {1}{2}( \gamma + \theta)}, \qquad f_{2(-)}(\theta) =
\frac {sinh
\frac
{1}{2}( \gamma -
\theta)}{sinh \frac {1}{2}( \gamma +
\theta)}
\end{equation}

 The reason for such a scope is that ( unlike
$p=\bar p$ for $N=(2p-1)$ )
for even $N$
there is $\it no$ index $i =\bar
i$. The successive stages of of the
construction of
of solutions in App.$A$
make it amply explicit how the presence of a $p
(=\bar p)$, along
with the
structure of the projectors in our nested sequence, constrains
the
coefficients to be simply exponentials. The generalization for $N= 2n
\quad
(n > 1)$ of
the hyperbolic and elliptic solutions displayed above
will be explored
elsewhere
implementing our basis of
projectors.

\section{Discussion:}
\setcounter{equation}{0}

In Ref.$1$ braid matrices were studied
systematically via their spectral
resolutions on
appropriate bases of
projectors. Such a study was already initiated in
previous works
( Ref.$2$
and Ref.$7$ ) and led to canonical factorization and
diagonalization in
Ref.$1$.
In Sec.$8$ of Ref.$1$ this approach was taken to its limit. In the
other
sections almost
all ${\it known}$ braid matrices of interest were
studied via spectral
resolutions. In
Sec.$8$ a basis of projectors ( called
a "nested sequence" ) with
particularly simple,
attractive properties was
hopefully presented for constructing new classes
of braid
matrices in all
dimensions. In such a basis, satisfying $(2.3)$, one has
$N^2$
matrices,
each $N^2\times N^2$ and with only constant elements ( see
$(2.2)$ and
$(6.1)$ ). They can
be considered as the most simple and
symmetric generalizations of
projectors appearing in
the $6$-vertex and the
$8$-vertex models. But the central question was not
addressed in
Ref.$1$.
Can such a basis of projectors be dressed up with suitable
coefficients
to
provide a braid matrix satisfying $(2.5)$ ? While the number
of
coefficients increases as
$N^2$ the number of trilinear constraints on
them corresponding to the
products of $N^3
\times N^3$ matrices increases
much faster. Hence the question. In this
paper we we
present an affirmative
answer and explicit solutions for all ${\it odd}$
$N$. The even-$N$
case
will be studied elsewhere.

  Let us note some basic features of our
solutions in the context of the
formulation in
Ref.$1$. The canonicaly
factorizable form of the coefficients $[1]$ give
\begin{equation}
\hat R
(\theta) = \sum_{i}\frac{f_i (\theta)}{f_i (-\theta)}
P_i
\end{equation}
This is evidently compatible with $(2.7)$ since
$$ e^{m
\theta} = (e^{\frac {1}{2}m \theta})(e^{-\frac {1}{2}m \theta})^{-1}$$

But in  Ref.$1$ we systematically extracted ( see the relevant
discussion
in Ref.$1$ )
the standard ( non-Baxterized ) braid matrices
satisfying
\begin{equation}
\hat{R}_{12}\hat{R}_{23}\hat{R}_{12}=
\hat{R}_{23}\hat{R}_{12}\hat{R}_{23}
\end{equation}
as the
limits
\begin{equation}
 \lim_{\theta \rightarrow \pm \infty} \hat R
(\theta) = (\hat R)^{\pm 1}
\end{equation}

 For our present class of
solutions however each coefficient $
e^{m_{\alpha}\theta}$
either diverges
or vanishes in the above limits. So rather than the
Baxterization of
a
preexisting $(7.2)$ to $(2.5)$ this class can be considered ( see Sec.$5$
)
to be an
exponentiation of
\begin{equation}
[[H_{12},H_{23}],H_{12}] =
[[H_{23},H_{12}],H_{23}]
\end{equation}

to
\begin{equation}
\hat{R}_{12}(\theta)\hat{R}_{23}(\theta+{\theta}')\hat{R}_{1
2}({\theta}')=
 \hat{R}_{23}({\theta}')\hat{R}_{12}({\theta+\theta}')
\hat{R}_{23}(\theta)
\end{equation}

since, as shown in Sec.$5$, the
passage

$$H \equiv \sum_{\alpha'} m_{\alpha'} P_{\alpha'} \quad
\rightarrow \hat
R(\theta) =
e^{\theta H}$$

correspopnds to one from
$(7.4)$ to $(7.5)$.

 One may compare this with the well-known so called
"classical" $r$- matrix
equation
obtained by expanding the $q$-dependent YB
matrix $R(\theta) (=P
\hat
R(\theta))$
satisfying

\begin{equation}
R_{12}(\theta)R_{13}(\theta+{\theta}')R_{23}({\theta}')=
 R_{23}({\theta}')R_{13}({\theta+\theta}')
R_{12}(\theta)
\end{equation}

in powers of $h(=ln q)$. One obtains for

$$
R_{q}(\theta) =I + 2hr(\theta) + O (h^2) $$

\begin{equation}
[r_{12}(\theta),( r_{13}(\theta+{\theta}')
+r_{23}({\theta}') )] +
[r_{13}({\theta+\theta}'), r_{23}({\theta}')] =
0
\end{equation}

 This has only single commutators. In our case there is
no $q$. Expanding
in powers of
$\theta$ we obtain as the first nontrivial
relation the equation $(7.4)$
with double
commutators and with the two
sides still directly related through the
interchange $(12)
\leftrightarrow
(23)$. In the extensive literature concerning $r$-matrices
one may note
in
particular a classification of solutions ( Ref.$8$ ). Our projectors
lead
to a solution of
$(7.4)$ with $\frac{1}{2}(N+3)(N-1))$ parameters for
$N=(2p -1)$. A more
general study,
starting from $(7.4)$ should be
worthwhile.

 We repeat a feature noted in Sec.$5$. Our class of solutions
has many
particularly simple aspects. But the number of projectors $(
P_{\alpha})$
and that of the
parameters $( m_{\alpha})$ increase as $N^2$
with the dimension. The degree
of the minimal
polynimial equation satisfied
by $\hat R(\theta)$ increases with them. This
is in sharp
contrast with
well-known cases corresponding to $q$-deformed  unitary,
orthogonal
and
symplectic cases. There the structures of the projectors are much
less
simple. But their
number does not increase with the dimension. As
noted below $(2.14)$, the
degree of the
minimal polynomial can be lowered
by allowing some of the free parameters
to coincide,
giving simpler
subcases. But our solution is more general.

 For $m^{+}_{ab} >m^{-}_{ab}$
all the nonzero elements of our $\hat
R(\theta)$ are
positive and hence can
be consistently interpreted as Boltzmann weights of
a
multistate
statistical model. In Sec.$11$ of Ref.$1$ the possibility of a
class of
multistate model
was briefly indicated and compared with one
proposed in Ref.$9$. ( See also
Sec.$4$ of
Ref.$10$. ) In both cases $(2N^2
- N)$ elements out of $N^4$ ones of $\hat
R(\theta)$ are
nonzero. Here we
have $(2N^2 -1)$ nonzero weights. Moreover the explicit
solution of
Ref.$9$
( and Ref.$10$ ) restricts the number of parameters as in the
$6$-vertex
model (
Sec.$6$ ). For our present class there is more scope in this
respect.

\bigskip

  It is a pleasure to thank Daniel Arnaudon. Using
a program, he verified
for the first
member of our hierarchy of solutions
that the constraints obtained here are
not only
sufficient  but also
necessary. This was reassuring.

\bigskip

\section {APPENDIX A. Solving
the braid
equation:}
\setcounter{equation}{0}
\renewcommand{\theequation}{A.\arabic{equation}}

 In $(2.6)$,
namely,

$$(\hat{R}(\theta))_{al,cm}(\hat{R}(\theta+{\theta}'))_{mn,ef}
(\hat{R}({\theta}'))_{lb,nd}$$
\begin{equation}
=(\hat{R}({\theta}'))_{cl,em}(\hat{R}({\theta+\theta}'))_{ab,ln}
(\hat{R}(\theta))_{nd,mf}
\end{equation}
corresponding to the site $(ab)\otimes (cd)\otimes (ef)$ one has to
implement the content of the ansatz $(2.4)$. From
$(2.2)$ and $(2.4)$ one obtains the following nonzero  elements of $\hat
R(\theta)$. The arguments $\theta$ is suppressed in $(A.2)$ to simplify
the notation and the subscripts correspond to the sites
$(ab)\otimes (cd)$.

$$ \hat R_{pp,pp}= 1$$
$$ \hat R_{pp,ii}= \frac{1}{2} ( f_{pi}^{(+)}+f_{pi}^{(-)} )
=\hat R_{pp,\bar{i} \bar{i}}$$
$$ \hat R_{pp,i\bar{i}}=
\frac{1}{2} ( f_{pi}^{(+)} -f_{pi}^{(-)} )
=\hat R_{pp,\bar{i} i}$$
$$ \hat R_{ii,pp}= \frac{1}{2} ( f_{ip}^{(+)} +f_{ip}^{(-)} )
=\hat R_{\bar{i}\bar{i},pp}$$
$$ \hat R_{i\bar{i},pp}= \frac{1}{2} ( f_{ip}^{(+)}
-f_{ip}^{(-)} )=\hat R_{\bar{i} i,pp}$$
$$ \hat R_{ii,jj}= \frac{1}{2} (f_{ij}^{(+)} +f_{ij}^{(-)} )
=\hat R_{\bar{i} \bar{i},\bar{j} \bar{j}}$$
$$\hat R_{i\bar{i},j \bar{j}}= \frac{1}{2} ( f_{ij}^{(+)} -f_{ij}^{(-)})
=\hat R_{ \bar{i}i,\bar{j}j}$$
$$ \hat R_{ii,\bar{j} \bar{j}}=\frac{1}{2} ( f_{i\bar{j}}^{(+)}
+f_{i\bar{j}}^{(-)} )=\hat R_{\bar{i}\bar{i},jj}$$
\begin{equation}
\hat R_{i\bar{i},\bar{j}j}= \frac{1}{2} (
f_{i\bar{j}}^{(+)}
-f_{i\bar{j}}^{(-)} )
=\hat
R_{\bar{i}i,j\bar{j}}
\end{equation}

These are the ${\it only }$ nonzero
elements, the total number being
$$  1+8(p-1)+8(p-1)^2 \quad = 2(2p -1)^2-1  \quad = 2N^2 - 1 $$

Note the following points:

$\bullet$   The
elements above all being situated on the diagonal and
the
antidiagonal
there are none of the type $\hat R_{i\bar{i},jj},\hat
R_{ii,j\bar{j}}$ and
so on.

$\bullet$  In the porduct $(ab)\otimes
(cd)\otimes (ef)$ for a given $a$,
$b$ can only be
$a$ or  $\bar{a}$ for
the coefficient to be nonzero. This holds also for
the other
pairs.

$\bullet$  Among $(a,b,c,d,e,f)$ the number of with  ( or without )
bar
must be even for
the coefficient to be nonzero. This is one consequence
of $(A.2)$. However,
in such
countings one must keep in mind that
$p=\bar{p}$.

 The preceding considerations simplify considerably the
computations as we
analyse
systematically the different classes of
$(ab)\otimes (cd) \otimes (ef)$
with novanishing
coefficients, lowering the
multiplicity of $(pp)$ in the triple product
above by steps.

 Case $(1)$:
The case $(pp)\otimes (pp) \otimes (pp)$ is trivial since
(A.$1$) reduces
to
        $$1=1$$
Case $(2)$: Next consider the classes ( with $(ab) \neq
(pp) $)
 $$ (1): (pp)\otimes (pp) \otimes (ab),$$
 $$ (2): (ab)\otimes (pp)
\otimes(pp),$$
 $$ (3): (pp)\otimes (ab) \otimes (pp)$$
 From our previous
remarks it follows that it is sufficient to consider the
possibilities
$$
(ab) = (ii), (i\bar{i})$$
Note also that in (A.$2$) $ \hat R _{pp,ii} =
\hat R _{pp,\bar i \bar i}$
and so on.

 For $(1)$, $(A.1)$ is easily seen
to reduce to

\begin{equation}
\hat R(\theta +\theta'))_{pp,ab} = \hat
R(\theta'))_{pp,ac}\hat
R(\theta))_{pp,cb}
\end{equation}

 Analogous
treatments of the subcases $(1),(2),(3)$ lead respectively
(
implementing
$(A.2)$ with $\epsilon =\pm$ and also both possibilities for
$(ab)$
mentioned above) to
the
constraints

\begin{equation}
f^{(\epsilon)}_{pi}(\theta +\theta') =
f^{(\epsilon)}_{pi}(\theta
)f^{(\epsilon)}_{pi}(\theta')
\end{equation}
\begin{equation}
f^{(\epsilon)}_{ip}(\theta +\theta') =
f^{(\epsilon)}_{ip}(\theta
)f^{(\epsilon)}_{ip}(\theta')
\end{equation}
$$f^{(+)}_{pi}(\theta )f^{(+)}_{pi}(\theta')f^{(+)}_{ip}(\theta +\theta')
+f^{(-)}_{pi}(\theta) f^{(-)}_{pi}(\theta')f^{(-)}_{ip}(\theta+\theta')$$
\begin{equation}
=
f^{(+)}_{pi}(\theta+\theta')f^{(+)}_{ip}(\theta )f^{(+)}_{ip}(\theta' )
+
f^{(-)}_{pi}(\theta +\theta')f^{(-)}_{ip}(\theta) f^{(-)}_{ip}(\theta')
\end{equation}

On implementing (A.$4$) and (A.$5$) one reduces (A.$6$)
to an identity.
Then from the
first two one obtains the
solutions
\begin{equation}
f^{(\epsilon)}_{pi}(\theta ) =
e^{m^{(\epsilon)}_{pi}\theta
}
\end{equation}
\begin{equation}
f^{(\epsilon)}_{ip}(\theta ) =
e^{m^{(\epsilon)}_{ip}\theta }
\end{equation}

the indeterminates
$m^{(\epsilon)}_{pi},m^{(\epsilon)}_{ip}$ being
independent parameters.

Continuing to reduce the multiplicity of $(pp)$ and remembering
the
restrictions implied
by  $(A.2)$  we start by considering successively
the cases
$$ (4): (pp)\otimes (ii) \otimes (jj)$$
$$ (5): (pp)\otimes
(i\bar {i}) \otimes (j \bar {j})$$
$$ (6):(pp)\otimes (i\bar {i}) \otimes
(jj)$$
The last one survives with nonzero coefficient since $p=\bar p$. We
present
directly the
the results, the derivations being straightforward.

Defining
$$ A_{ab}(\theta) \equiv f^{(+)}_{ab}(\theta )
+f^{(-)}_{ab}(\theta ) , \qquad
B_{ab}(\theta) \equiv f^{(+)}_{ab}(\theta )
-f^{(-)}_{ab}(\theta ) $$

one obtains respectively from the above
cases
\begin{equation}
A_{ij} (\theta +\theta') = f^{(+)}_{ij}(\theta
)f^{(+)}_{ij}(\theta' ) +
f^{(-)}_{ij}(\theta )f^{(-)}_{ij}(\theta'
)
\end{equation}
$$A_{pi}(\theta) A_{pi}(\theta') B_{ij}(\theta +\theta')
+B_{pi}(\theta)
B_{pi}(\theta')
B_{i\bar{j}}(\theta +\theta')
$$
\begin{equation}
= A_{pi}(\theta +\theta') (B_{ij}(\theta)
A_{ij}(\theta') +A_{ij}(\theta)
B_{ij}(\theta')
)
\end{equation}
$$A_{pi}(\theta) B_{pi}(\theta') A_{ij}(\theta +\theta')
+B_{pi}(\theta)
A_{pi}(\theta')
A_{i\bar{j}}(\theta +\theta')
$$
\begin{equation}
= B_{pi}(\theta +\theta') (A_{ij}(\theta')
A_{i\bar{j}}(\theta)
+B_{ij}(\theta')
B_{i\bar{j}}(\theta)
)
\end{equation}

  Takingt account of$(A.4)$ and $(A.5)$ ( and hence of
$(A.7)$ and $(A.8)$
) and noting
that keeping $(\theta +\theta')$ fixed one
can vary $\psi$ in
$$ \theta = \phi +\psi, \quad \theta' = \phi - \psi
$$
one finds that the last three equations are satisfied
if
\begin{equation}
f_{i\bar{j}}^{(\epsilon)}(\theta) =
f_{ij}^{(\epsilon)}(\theta)
\end{equation}
and
\begin{equation}
f_{ij}^{(\epsilon)}(\theta)f_{ij}^{(\epsilon)}(\theta') =
f_{ij}^{(\epsilon)}(\theta+\theta')
\end{equation}

These are found to be ${\it necessary}$ and
${\it sufficient}$.
Hence
\begin{equation}
f_{i\bar{j}}^{(\epsilon)}(\theta) =
f_{ij}^{(\epsilon)}(\theta) =
e^{m_{ij}^{(\epsilon)}\theta}
\end{equation}

Permutation of the factors of the cases $(4,5,6)$ above ( such
as
$(ii)\otimes (pp)
\otimes (jj)$ and so on) can be shown to lead to no
supplementary constraints.

 Finally one considers the cases
$$ (ab)\otimes
(cd) \otimes (ef)$$
where no factor is $(pp)$. For each subcase the
constraints implied by
(A.$1$) along with
(A.$2$) are easily extracted. It
is found that they are ${\it all}$
satisfied by
implementing $(A.12)$ and
$(A.13)$. Since the subcases are treated quite
similarly, it is
sufficient
to display two of them. We present again only the the final
steps. For
$$
(ii)\otimes (jj) \otimes (kk)$$
with no barred index, (A.$13$) reduces
(A.$1$), in terms of $A_{ab}$
defined above, to
\begin{equation}
L.H.S. =
\frac{1}{4}A_{ij}(\theta +\theta')A_{jk}(\theta +\theta') =
R.H.S.
\end{equation}
Similarly, for
$$ (i\bar{i})\otimes (jj) \otimes
(k\bar{k})$$
one obtains finally
\begin{equation}
L.H.S. =
\frac{1}{4}B_{ij}(\theta +\theta')B_{jk}(\theta +\theta') =
R.H.S.
\end{equation}

 In both cases, apart from the exponential form for
each $f$, $(A.14)$ is
essential.
Thus we have verified the solution
announced in $(2.7)$ and $(2.8)$. It is
instructive to
compute explicitly
the case $(2.9)$ where one has only $(i,\bar i,p)$ with
$$ i=1, \quad
p=2$$
One finds that $(2.10)$ is sufficient. Moreover, if one
sets
\begin{equation}
 m_{1 \bar 1}^{(\epsilon)} \neq  m_{1
1}^{(\epsilon)}
\end{equation}
so that $a_{\pm}$ is not repeated as in
$(2.9)$, the braid equation is
${\it not}$
satisfied. This is an example of
the necessity of $(A.14)$.

 As a check, the solution for $N=3$ was also
obtained ( instead of directly
using (A.$1$)
and (A.$2$) ) by computing the
triple tensor products of the projectors in
$(2.5)$.

\section
{APPENDIX B. L-operators and transfer matrices (
fundamental
representations
)
:}

\setcounter{equation}{0}
\renewcommand{\theequation}{B
.\arabic{equation}}

 Here we collect together some known results ( citing
sources below )
coherently with our
notations and conventions and emphasize
certain aspects arising in the
presence of the
spectral parameter
$\theta$.

 For non-Baxterized braid matrices ( without $\theta$ )
satisfying
\begin{equation}
\hat{R}_{12}\hat{R}_{23}\hat{R}_{12}=
\hat{R}_{23}\hat{R}_{12}\hat{R}_{23}
\end{equation}

the FRT equations for
the $L$-operators ( eqn.$(2.3)$ of Ref.$11$ ) can be
expressed
in our
notations as
\begin{equation}
\hat R L_{2}^{\pm} L_{1}^{\pm} = L_{2}^{\pm}
L_{1}^{\pm}\hat R
\end{equation}
\begin{equation}
\hat R L_{2}^{+}
L_{1}^{-} = L_{2}^{-} L_{1}^{+}\hat R
\end{equation}

  Here $\hat R$ is a
$N^2 \times N^2$ matrix for any $N$ and
$$L_1 = L\otimes I_{N\times N},
\quad L_2 = I_{N\times N} \otimes L$$

 Writing these in terms of
components ( as will be done below for the
$\theta$-dependent
case ) it can
be shown that the lowest dimensional realizations of the
$N^2$
blocks
$L^{\pm}_{ab}$ ( each $N\times N$ ) can be obtained in our notations
i,e, with
\begin{equation}
\hat R = \hat R_{ab,cd} (ab)\otimes
(cd)
\end{equation}
as
\begin{equation}
(L^{+}_{ab})_{cd} = \hat R
_{ad,cb}
\end{equation}
\begin{equation}
(L^{-}_{ab})_{cd} = \hat R^{-1}
_{ad,cb}
\end{equation}
or
\begin{equation}
L^+ = \hat R P = (PR)P =
R_{21}
\end{equation}
\begin{equation}
L^-={\hat R }^{-1}P
=(R^{-1}P)P=R^{-1}
\end{equation}

 Apart from differences of notations and
conventions these correspond ( to
cite only one
source ) to eqns.$(4.9)$ of
Ref.$12$. In the familiar $L^{\pm}$ of $Sl_{q}(2)$,
implementing $2\times
2$ realizations of $(q^{\pm H}, X_{\pm})$ one obtains
$(B.7)$ and
$(B.8)$,
which however hold for ${\it any}$ $\hat R$ satisfying $(B.1)$.

 Now let
us introduce $\theta$. Corresponding to $(B.1)$ and $(B.2)$ one
now
has
respectively
\begin{equation}
\hat{R}_{12}(\theta -
\theta')\hat{R}_{23}(\theta)\hat{R}_{12}(\theta')=
\hat{R}_{23}(\theta')\hat{R}_{12}(\theta) \hat{R}_{23}(\theta
-\theta')
\end{equation}
\begin{equation}
\hat{R}(\theta -
\theta')L_{2}^{\pm}(\theta)L_{1}^{\pm}(\theta')=
L_{2}^{\pm}(\theta')L_{1}^{\pm}(\theta) \hat{R}(\theta-\theta')
\end{equation}
( The corresponding situation for (B.$3$) will be
discussed below. )

In terms of components one
writes
$$(\hat{R}(\theta
-\theta'))_{al,cm}(\hat{R}(\theta))_{mn,ef}(\hat{R}({\theta}'))_{lb,nd}$$
\begin{equation}
=(\hat{R}({\theta}'))_{cl,em}(\hat{R}(\theta))_{ab,ln}
(\hat{R}(\theta
-\theta'))_{nd,mf}
\end{equation}
$$(\hat{R}(\theta-\theta'))_{a
l,cm}(L^{\pm}(\theta)_{mf})_{en}
(L^{\pm}({\theta}')_{ld})_{nb}$$
\begin{equation}
=(L^{\pm}({\theta}')_{cm})_{el}(L^{\pm}(\theta)_{an})_{lb}
(\hat{R}(\theta-\theta'))_{nd,mf}
\end{equation}

 One finds that ( considering
$L^{+}(\theta)$ to start with
)
\begin{equation}
(L^{+}(\theta)_{ab})_{cd}=(\hat{R}(\theta))_{ad,cb}
\end{equation}
or
\begin{equation}
L^{+}(\theta)=\hat{R}(\theta)P \quad =
PR(\theta)P \quad = R_{21}(\theta)
\end{equation}
is a solution. This is
strictly analogous to (B.$7$). The same solution
evidently holds
for
$L^{-}(\theta)$. But if one wants to avoid the
degeneracy
\begin{equation}
L^{-}(\theta)=L^{+}(\theta)
\end{equation}
can
one obtain a different solution for $L^{-}(\theta)$ analogous to
(B.$8$) ?
We show
below in a particularly transparent fashion that there is an
obstruction if
one directly
generalizes (B.$3$) as
\begin{equation}
\hat R
(\theta -\theta') L_{2}^{+}(\theta) L_{1}^{-}(\theta')
=
L_{2}^{-}(\theta')
L_{1}^{+}(\theta)\hat R(\theta -
\theta')
\end{equation}
When this is further generalized by introducing a
central operator in the
argument of
$\hat R (\theta -\theta')$ on one side
( or in a different fashion on each
side ) and thus
distinguish the two
arguments, there can be a way out. ( Ref.$13$ is a
review article
citing
numerous sources. Particularly relevant is Sec.$2.1.4$.) But let
us
consider the
consequences of (B.$16$) combined with (B.$13$), (B.$14$)
and the basic
properties
$(2.15)$, i,e,
\begin{equation}
 \hat R(-\theta)
={\hat R}^{-1}(\theta), \qquad  \hat R(0) =I
\end{equation}
>From (B.$14$)
and (B.$17$),
\begin{equation}
L^{+}(0) = P, \qquad L^{+}(0)_{ab} =
(ba)
\end{equation}

(This has no counterpart for (B.$7$).)

Hence setting
$\theta =0$ in (B.$16$), using (B.$16$) and (B.$17$) and
then
writing
$\theta$ for $\theta'$ one obtains
\begin{equation}
{\hat
R}^{-1} (\theta) P_{2} L_{1}^{-}(\theta) = L_{2}^{-}(\theta)
P_{1}{\hat
R}^{-1}(\theta)
\end{equation}

Writing (B.$19$) in terms of components
analogously to (B.$12$) one
obtains
\begin{equation}
(\hat{R}(-\theta))_{al,cm}(\delta_{mn}
\delta_{fe})(L^{-}({\theta})_{ld})_{nb}=
(L^{-}({\theta})_{cm})_{el}(\delta_{ab}
\delta_{nl})(\hat{R}(-\theta))_{nd,mf}
\end{equation}
or
\begin{equation}
({\hat
{R}}^{-1}(\theta))_{al,cm}(L^{-}({\theta})_{ld})_{mb}\delta_{ef}=
\delta_{ab}(L^{-}({\theta})_{cm})_{el}({\hat{R}}^{-1}(\theta))_{ld,mf}
\end{equation}
Hence, finally,
\begin{equation}
({\hat R}^{-1}(\theta)
L^{-}(\theta)P)\otimes I = I\otimes (L^{-}(\theta)P{\hat R}^{-1}(\theta))
\end{equation}

For
\begin{equation}
L^{-}(\theta) =L^{+}(\theta) = \hat
R(\theta)P, \qquad L^{-}(\theta)P=\hat R(\theta)
\end{equation}
(B.$22$) is
trivially satisfied ( furnishing a convincing check ). But a
distinct solution
for $L^{-}(\theta)$ reducing ( say, as $\theta \rightarrow \infty$ ) to
(B.$8$) is no longer available in the general case if (B.$16$) is
strictly maintained. One obtains (B.$8$) easily from the symmetry of
(B.$1$) under inversion since,unlike for (B.$16$), the orders of
$(\theta,\theta')$ on each side do not enter in that context. But even
apart from that (B.$18$) now imposes the constraint (B.$22$), linear
in $L^{-}(\theta)$. We do not
consider in this paper generalizations of (B.$16$) leading to quasi-Hopf
structures for consistent coproducts.

We
now consider transfer matrices and note how the lowest dimensional
representations can be extracted from those of the
$L$-operators. The transfer matrix $t(\theta)$ has to satisfy
\begin{equation}
\hat R(\theta -\theta') (t(\theta)
\dot{\otimes} t(\theta')) =(t(\theta')
\dot{\otimes}
t(\theta))\hat
R(\theta -\theta')
\end{equation}
where $\dot \otimes$, combining tensor
and matrix products leads to
\begin{equation}
(t(\theta) \dot{\otimes}
t(\theta')) =(t(\theta) \otimes I)(I \otimes
t(\theta'))
\quad
=t_{1}(\theta). t_{2}(\theta')
\end{equation}
Writing (B.$24$) as
$$(P\hat R(\theta -\theta')P)(P (t_{1}(\theta)P)(Pt_{2}(\theta')P)$$
\begin{equation}
=(Pt_{1}(\theta')P)(Pt_{2}(\theta)P)(P\hat R(\theta
-\theta')P)
\end{equation}
or
\begin{equation}
{\hat R}_{21}(\theta
-\theta') t_{2}(\theta) t_{1}(\theta') =t_{2}(\theta')
 t_{1}(\theta){\hat
R}_{21}(\theta -\theta')
\end{equation}
Now comparing (B.$27$) with
(B.$10$) and (B.$14$) one finds the solution
\begin{equation}
t(\theta) =
(P\hat R(\theta) P)P \quad = P\hat R(\theta) \quad =
R(\theta)
\end{equation}
In absence of $\theta$, i,e, for (B.$1$), this
corresponds ( with some
notational
differences ) to the realization
$\rho_{+}$ of of eqn.$(4.5)$ of Ref.$12$. But
corresponding to (B.$8$),
unavailable in our context, there is anther
realization
$\rho_{-}$ in
Ref.$12$. We are concerned only with (B.$28$). Products
analogous to
$(4.8)$
of our Sec.$4$ lead to higher dimensional transfer matrices
corresponding
to longer chains
as successive sites are
added.

\bigskip

\bigskip

\bibliographystyle{amsplain}

\end{document}